\documentclass[11pt]{article}
\usepackage[margin=1.0in]{geometry}
\usepackage{amsthm}
\usepackage{amssymb,amsfonts,amsmath,latexsym}
\usepackage{cite}
\usepackage{graphicx}
\usepackage{hyperref}

\newcounter{count}
\theoremstyle{plain}
\newtheorem{definition}[count]{Definition}
\newtheorem{theorem}[count]{Theorem}
\numberwithin{count}{section}

\newcommand{\R}{\mathbb{R}}
\newcommand{\RP}{\mathbb{RP}}
\newcommand{\CP}{\mathbb{CP}}
\newcommand{\C}{\mathbb{C}}

\renewcommand{\L}{\Lambda}

\newcommand{\Cov}{\operatorname{Cov}}
\newcommand{\Var}{\operatorname{Var}}

\newcommand{\ra}{\rightarrow}

\title{A Schur transform for spatial stochastic processes}
\author{James Mathews\footnote{Memorial Sloan-Kettering Cancer Center (\href{}{mathewj2@mskcc.org})}}
\begin{document}

\maketitle

\begin{abstract} The variance, higher order moments, covariance, and joint moments or cumulants are shown to be special cases of a certain tensor in $V^{\otimes n}$ defined in terms of a collection $X_1,...,X_n$ of $V$-valued random variables, for an appropriate finite-dimensional real vector space $V$. A statistical transform is proposed from such collections--finite spatial stochastic processes--to numerical tuples using the Schur-Weyl decomposition of $V^{\otimes n}$. It is analogous to the Fourier transform, replacing the periodicity group $\mathbb{Z}$, $\mathbb{R}$, or $U(1)$ with the permutation group $S_{n}$. As a test case, we apply the transform to one of the datasets used for benchmarking the Continuous Registration Challenge, the thoracic 4D Computed Tomography (CT) scans from the M.D. Anderson Cancer Center available for download from DIR-Lab. Further applications to morphometry and statistical shape analysis are suggested.
\end{abstract}



\tableofcontents

\section{Introduction}

In engineering and the sciences there are many situations in which one measures the locations of $N$ points which are matched across a group of $n$ different objects. For example, anatomical landmarks in multi-atlas medical images, skeletal landmarks of specimens in the fossil record, or a time series of several point-motion trackers on a moving subject. What are the natural features of such data? Such features could be used to train predictive models, to feed back into control systems, or to simply re-present the data so that an investigator can discern patterns.

The geometry and dynamics of such data are normally quantified by measures like volumetric deformation, deformation tensor eigenvalues, kinetic energy, total displacement along paths of travel, path curvature, and the summary statistics of such quantities, as well as more inherently statistical measures like higher-order joint moments.

This article proposes a class of features arising from classical tensor algebra and representation theory that are natural generalizations of elementary statistics, like the covariance or joint moments or cumulants, and that are substantially coordinate-independent. We interpret the given data in the context of probability theory by positing that the locations of the $N$ points on $n$ objects are samples from the joint distribution of $n$ different random spatial or vector variables valued in a vector space $V$. The expected value of the tensor product of these variables is a tensor in the classical sense, an element of the tensor product space $V^{\otimes n}=V\otimes \cdots \otimes V$. It is at this point that the multivariate case shows itself to be much richer and much more challenging than the univariate case: The tensor product of 1-dimensional vector spaces is always another 1-dimensional vector space, so that tensor language can be entirely avoided in the theory of univariate statistics, while the $n$-fold tensor product of a $k$-dimensional vector space for $k>1$ is literally exponentially larger, of dimension $k^{n}$.

There is an additive decomposition of elements of $V^{\otimes n}$ called the \emph{Schur-Weyl decomposition} (not to be confused with the multiplicative Schur decomposition of a matrix). The Schur-Weyl decomposition arises from the consideration of the possible symmetry types with respect to permutations of the $n$ factors. Or, what turns out to be equivalent, the possible symmetry types with respect to the linear transformations of $V$, the elements of the general linear group $GL(V)$. These symmetry types in some sense interpolate between the fully-symmetric type (as in symmetric matrices) and fully-antisymmetric type (as in anti-symmetric matrices). The new features that we emphasize are the norms of the components of this decomposition, which we call Schur amplitudes, thought of as invariants for a restricted, orthogonal group. One could also consider any other invariants of these components.


We will recall the part of the classical representation theory of $GL(V)$ and the permutation group $S_n$ that is required for the actual numerical calculation of the proposed features. As it turns out, the formula for the components of the Schur-Weyl decomposition is a convolution of exactly the same sort as the Fourier transform. The main difference is that the convolution kernel for the Fourier transform, interpreted as the character functions for the irreducible representations of the unitary group $U(1)$ or the additive groups $\R$ or $\mathbb{Z}$, is replaced with a convolution kernel given by the character functions $\chi_{\lambda}$ for the irreducible representations $\lambda$ of the symmetric group $S_n$, also known as the permutation group on $n$ symbols.

After explaining the algorithm in detail, we present the results of the transform applied to 4DCT lung motion data as a test case. This very preliminary example is used to give possible interpretations of the Schur transform in applications.

\section{Theory}

Let $V$ be a finite-dimensional real vector space and let $X_1,...,X_n$ be $V$-valued random variables on a fixed probability space $\Omega$. We assume that the $X_i$ are jointly continuous with joint probability density $d(v_1,...,v_n)$ on $V\times \cdots \times V$. We use the notation $S^{l}V$ for the $l^{th}$ symmetric power of $V$.

The following definition is motivated by the need for a uniform treatment of the types of quantities that will be mentioned in Table \ref{specialcases}.

\begin{definition} \label{covt}The \emph{covariance tensor}, denoted $t$, is
\[\Cov(X_1,...,X_n):=\int_{V\times\cdots\times V} (v_1-\mu_1) \otimes \cdots \otimes (v_n-\mu_n) \quad d(v_1,...,v_n) \in V^{\otimes n}\]

where $\mu_i:=\mathbb{E}(X_i)\in V$.

The \emph{covariance tensor of type $(l_i)_{i=1}^{n}$}, denoted $t_{(\l_i)}$, where each $l_i$ is a positive integer, is
\begin{align*}
&\Cov_{(l_1,...,l_n)}(X_1,...,X_n):=\\
&\quad\Cov(\underbrace{X_1,...,X_1}_{l_1},...,\underbrace{X_n,...,X_n}_{l_n})\in S^{l_1}V\otimes \cdots \otimes S^{l_n}V\subset V^{\otimes \sum_i^{n} l_i}
\end{align*}
\end{definition}

The covariance tensor is still defined, by the same formula, if each $X_i$ is valued in an affine space $A$ for the vector space $V$ rather than $V$ itself. This is because in this case, each $v-\mu_i$ is still defined as an element of $V$ since $v,\mu_i \in A$. Also, the definition of the covariance tensor does not require a choice of inner product, a basis, or any other additional structure on $V$.

In some situations one has meaningful origin or reference points for each variable $X_i$ which are not necessarily the means $\mu_i$. In this case it might be desirable to replace the means $\mu_i$ in Definition \ref{covt} with these reference points. The resulting tensor may be called \emph{non-central}, by analogy with the non-central ordinary moments.

Several special cases of the covariance tensor $t$ are summarized in Table \ref{specialcases} and described as follows.

If $V=\R$, there is a canonical isomorphism $V^{\otimes n}\cong \R$. If in addition $n=2$, $t$ is the ordinary univariate covariance $\Cov(X_1,X_2)$. On the other hand, for $n=1$, $t_{(2)}$ is the variance $\Var(X_1)$. $t_{(3)}$, $t_{(4)}$, etc. are the higher moments of $X_1$. If $n>2$, $t$ is the joint moment of the $X_i$.

If $\dim V > 1$ and $n=2$, the matrix of $t$ with respect to a basis for $V$ is the covariance matrix of $X_1$ and $X_2$, the matrix $\{\Cov((X_{1})_{\alpha},(X_{2})_{\beta})\}_{\alpha,\beta = 1 ... \dim V}$.

If $\dim V > 1$ and $n=1$, the matrix of $t_{(2)} \in S^{2}V$ with respect to a basis for $V$ is the variance matrix of a single vector-valued variable $X_1$. This $t_{(2)}$ is also called the \emph{inertia tensor} of $X_1$, especially when $\dim V=3$.

For general $\dim V$ we shall call $t_{(3)}$, $t_{(4)}$, $t_{(5)}$, etc. the \emph{spatial higher moments} of $X_1$ in $S^{3}V$, $S^{4}V$, $S^{5}V$, etc. We have singled out the cases $t_{(4)}\in S^{4}\R^{2}$ and $t_{(3)}\in S^{3}\R^{3}$ because in both cases there is a classical numerical invariant, known as the $J$-invariant, which captures essential features of the geometry of the tensor and which is well-defined under arbitrary affine transformations of $V=\R^{2}$ or $V=\R^{3}$ respectively. For an element $p\in S^{4}\R^{2}$ regarded as a homogeneous quartic polynomial in 2 variables, the $J$-invariant describes the set of the 4 complex roots of $p$ in $\CP^{1}$ up to projective transformations. For an element $q\in S^{3}\R^{3}$ regarded as a homogeneous cubic polynomial in 3 variables, the $J$-invariant describes the algebraic isomorphism type of the elliptic curve described by $q$ in the projective plane $\RP^{2*}$.

\begin{table}[h]
\centering
\begin{tabular}{| c | c | c | c | c |} \hline
                             & $\dim V$ & $n$ & type $(l_i)$& tensor type              \\ \hline
covariance                   &   1      & 2   & $(1,1)$     & $V \cong \R$             \\ \hline
variance                     &   1      & 1   & $(2)$       & $S^{2}V \cong \R$        \\ \hline
$k^{th}$ higher moment       &   1      & 1   & $(k)$       & $S^{k}V \cong \R$        \\ \hline
joint moment                 &   1      & $>2$& $(1,...,1)$ & $V^{\otimes n}\cong \R$  \\ \hline
covariance matrix            &   $>1$   & 2   & $(1,1)$     & $V\otimes V$             \\ \hline
variance matrix              &   $>1$   & 1   & $(2)$       & $S^{2}V$                 \\ \hline
inertia tensor               &   3      & 1   & $(2)$       & $S^{2}V\cong S^{2}\R^{3}$\\ \hline
$4^{th}$-order moment in 2d  &   2      & 1   & $(4)$       & $S^{4}V\cong S^{4}\R^{2}$\\ \hline
$3^{rd}$-order moment in 3d  &   3      & 1   & $(3)$       & $S^{3}V\cong S^{3}\R^{3}$\\ \hline
\end{tabular}
\caption{\label{specialcases}Special cases of the covariance tensor.}
\end{table}

In general $t$ and $t_{(l_i)}$ are not expressible in terms of matrix algebra operations like sum, product, transpose, inverse, and trace, and they are not among the familiar objects from probability theory or statistics (see e.g. \cite{wass} or \cite{tuck}). Certain cases of this general covariance tensor are discussed in \cite{tensormethodsstat} in the fully tensorial context, but no systematic description by means of symmetries is attempted.

\vspace{1pc}

Let $\lambda$ denote an arbitrary positive integer partition of the integer $n$ (e.g. $\lambda=(4,3,1)$ for $n=8$). Recall that there is a one-to-one correspondence between the set of $\lambda$ (also described by Young diagrams with $n$ boxes or by conjugacy classes in $S_n$) and the isomorphism classes of irreducible representations of the permutation group $S_n$ (\cite{fultonharris} page 44). Thus $\lambda$ will also be used to denote such an isomorphism class, and $M_\lambda$ will denote an $S_n$ representation in the isomorphism class $\lambda$. $\chi_{\lambda}$ denotes the character function of $S_n$ associated with $M_{\lambda}$, that is $\chi_{\lambda}(\sigma)$ is the trace of the linear map $M_{\lambda}\ra M_{\lambda}$ by which the element $\sigma\in S_n$ acts. Similarly there is a one-to-one correspondence between the set of $\lambda$ and certain functors $V\mapsto S^\lambda(V)$, called Schur functors, where $S^\lambda(V)$ is a certain representation of $GL(V)$ (\cite{fultonharris} page 75).

Each element $\sigma\in S_n$ acts as a linear transformation of $V^{\otimes n}$ to permute the tensor factors. The behavior of the elements of $V^{\otimes n}$ under such transformations can be used to classify and decompose such elements. Frequently, a representation $W$ of a group $G$ has a unique decomposition as a direct sum of $G$-irreducible subrepresentations, but one is not always so lucky. As a representation of $S_n$, $V^{\otimes n}$ does not have a unique decomposition as a direct sum of irreducible subrepresentations as soon as $n>2$. The same is true with respect to the group $GL(V)$ rather than $S_n$. However, the following classical theorem provides an alternative. For this theorem only, we assume that $V$ is a finite-dimensional \emph{complex} vector space.

\begin{theorem}\label{schurweyl}
\leavevmode
\begin{enumerate}
	\itemsep0em
\item{\label{decomp}The irreducible decomposition of $V^{\otimes n}$ as a representation of $S_n\times GL(V)$ is unique, written $V^{\otimes n}= \bigoplus_\lambda I_\lambda$, where the $I_\lambda$ are the irreducible subrepresentations and $\lambda$ ranges over the positive integer partitions of $n$.}
\item{\label{comp}$I_\lambda$ is the sum of all $S_n$-subrepresentations of $V^{\otimes n}$ isomorphic to $M_\lambda$.}
\item{\label{comp2}$I_\lambda \cong M_\lambda \otimes S^\lambda(V)$ as $S_n\times GL(V)$ representations.}
\item{\label{projf}The projection $\pi_\lambda:V^{\otimes n}\rightarrow I_\lambda\subset V^{\otimes n}$ is
\[\pi_\lambda : t \mapsto \frac{\chi_{\lambda}(\operatorname{identity})}{n!}\sum_{\sigma \in S_n} \chi_{\lambda}(\sigma^{-1})\sigma(t)\]
where $\chi_{\lambda}$ denotes the character function of $M_\lambda$ and $\sigma(t)\in V^{\otimes n}$ denotes the action of $\sigma$ on $t$.}
\end{enumerate}
\end{theorem}

Theorem \ref{schurweyl} records just a small part of a large general theory, the part which is needed here. Theorem \ref{schurweyl} goes all the way back to Issai Schur's 1901 dissertation \cite{schur1901} and his later article \cite{schur1927}. Hermann Weyl's 1939 book \cite{weyl1939classical}, chapter IV, is a well-known exposition in English. See \cite{berget} page 81 for discussion and further references.

(\emph{Dimensions}). It is no trivial matter to establish the number of degrees of freedom for each component of the decomposition in Theorem \ref{schurweyl}, the dimensions of the $I_{\lambda}$. Nevertheless these dimensions are classically known. The dimension of $M_\lambda$ is provided by the hook-length formula in terms of the Young diagram of type $\lambda$. The dimension of $S^{\lambda}(V)$ is given by the value of character function for $S^{\lambda}(V)$ on the identity element of $GL(V)$, expressed in terms of the eigenvalues of an element $L\in GL(V)$ by the $\lambda^{th}$ Schur polynomial. Then $\dim I_\lambda = \dim M_{\lambda} \cdot \dim S^{\lambda}(V)$. 

Surprisingly, it follows from the Frobenius formula (\cite{fultonharris} page 49) that the values of the character functions $\chi_\lambda$ for $S_n$ are integers. This implies that the decomposition exists for real $V$, so we did not need to pass to the complexification $V_{\C}$ after all. It also implies that in principle the projection by $\pi_\lambda$ can be performed exactly by a computer, without special commutative ring calculations and without floating-point approximation.

\begin{table}[h]
\setlength{\tabcolsep}{0.02em}
\begin{center}
\scalebox{0.85}{
\begin{tabular}{ |c| c| c| c| c| c| c| c| c| c| c| c|}\hline
c.c. size& 1 & 15 & 45 & 15 & 40 & 120 & 40 & 90 & 90 & 144 & 120 \\ \hline
c.c. rep.& () & (12) & (12)(34) & (12)(34)(56) & (123) & (123)(45) & (123)(456) & (1234) & (1234)(56) & (12345) & (123456) \\ \hline
 & 1 & -1 & 1 & -1 & 1 & -1 & 1 & -1 & 1 & 1 & -1 \\  \cline{2-12}
 & 5 & -3 & 1 & 1 & 2 & 0 & -1 & -1 & -1 & 0 & 1 \\  \cline{2-12}
 & 9 & -3 & 1 & -3 & 0 & 0 & 0 & 1 & 1 & -1 & 0 \\  \cline{2-12}
 & 5 & -1 & 1 & 3 & -1 & -1 & 2 & 1 & -1 & 0 & 0 \\  \cline{2-12}
 & 10 & -2 & -2 & 2 & 1 & 1 & 1 & 0 & 0 & 0 & -1 \\  \cline{2-12}
 & 16 & 0 & 0 & 0 & -2 & 0 & -2 & 0 & 0 & 1 & 0 \\  \cline{2-12}
 & 5 & 1 & 1 & -3 & -1 & 1 & 2 & -1 & -1 & 0 & 0 \\  \cline{2-12}
 & 10 & 2 & -2 & -2 & 1 & -1 & 1 & 0 & 0 & 0 & 1 \\  \cline{2-12}
 & 9 & 3 & 1 & 3 & 0 & 0 & 0 & -1 & 1 & -1 & 0 \\  \cline{2-12}
 & 5 & 3 & 1 & -1 & 2 & 0 & -1 & 1 & -1 & 0 & -1 \\  \cline{2-12}
 & 1 & 1 & 1 & 1 & 1 & 1 & 1 & 1 & 1 & 1 & 1 \\ \hline
\end{tabular}}
\caption{The character table of the symmetric group $S_6$, provided for illustration. Computed with SageMath. Each column is: the size of a conjugacy class, a representative element of the conjugacy class, and the values taken by each character on the class. Each row corresponds to one character $\lambda$, one partition of 6. Note that the first column is the dimension of the vector space $M_{\lambda}$.}
\end{center}
\end{table}

The projection formula in Theorem \ref{schurweyl} part \ref{projf} is essentially the formula for the Fourier transform, as clarified in Table \ref{ftrans}.

\begin{table}[h]
\centering
\begin{tabular}{| l | l | l | l |} \hline
group        & dual group   & characters                           & Fourier transform \\ 
$G$          & $\widehat{G}$& $\chi_{\widehat{g}}:G\ra \C$         & $L^{2}(G)\ra L^{2}(\widehat{G})$\\ \hline
$U(1)$       & $\mathbb{Z}$ & $\chi_k(\theta)=e^{2\pi i\theta k}$  & $f(\theta) \mapsto (\int_{[\theta]\in U(1)}\chi_{k}(-\theta)f(\theta)d\theta)_{k\in \mathbb{Z}}$\\ \hline
$\mathbb{Z}$ & $U(1)$       & $\chi_{\theta}(k)=e^{2\pi i\theta k}$& $f(k) \mapsto \sum_{k\in \mathbb{Z}}\chi_{\theta}(-k)f(k)$\\ \hline
$\R$         & $\R$         & $\chi_{c'}(c)=e^{2\pi i c c'}$       & $f(c)\mapsto \int_{c\in \R}\chi_{c'}(-c)f(c)dc$\\ \hline
\end{tabular}
\caption{\label{ftrans}In all 3 standard cases of the Fourier transform, involving (1) the unitary group $U(1)=\{e^{2\pi i\theta}|\theta\in \R\}\subset \C$ (Fourier series from a periodic function), (2) the integers $\mathbb{Z}=\{k\}$ (a periodic function from Fourier series), and (3) the real numbers $\R=\{c\}$ (ordinary Fourier transform), the transform is given by the formula $f(g)\mapsto \mathcal{F}(f)(\widehat{g})=\int_{g\in G}\chi_{\widehat{g}}(-g)f(g)dg$.}
\end{table}

We are now ready to define our main object of study, the Schur transform and some of its variants.

\begin{definition} Let $t=\Cov(X_1,...,X_n)$, for $V$-valued random vectors $X_i$.

\begin{enumerate}
	\itemsep0em
\item{The \emph{Schur transform} of $(X_1,...,X_n)$ is the tuple of components in the Schur-Weyl decomposition of $t$, $(t_\lambda)=(\pi_{\lambda}(t))$.}
\item{Fix an inner product on $V$ and consider the induced inner product and norm on $V^{\otimes n}$. The \emph{Schur amplitudes} of $t$ are the norms $(|t_\lambda|)$.}
\item{Suppose that $m$ random vectors $(Y_1,...,Y_m)$ valued in $V$ are given, for $m>n$. The \emph{$n$-factor Schur content} is the set (distribution) of the Schur amplitudes of all $n$-element subsets of the $Y_i$.}
\item{Suppose that $m$ random vectors $(Y_1,...,Y_m)$ valued in $V$ are given, for $m>n$, as above. The \emph{sequential $n$-factor Schur content} is the set (distribution) of the Schur amplitudes of all consecutive $n$-element subsets of the $Y_i$.}
\end{enumerate}
\end{definition}

Notice that the Schur amplitudes are independent of the order of the variables $X_1,...,X_n$, so that the $n$-factor Schur content is well-defined. Also, if all $X_i$ are equal to some fixed $X$, the Schur transform has only one non-zero entry, which is equal to $\Cov_{(n)}(X)\in S^{n}V$.

Each component of the Schur transform has its own interpretation as an encoding of some aspect of the geometric arrangement of the $X_i$. For example, the $\L^{n}V$ component quantifies the tendency of the values of the $X_i$, as displacements from the means of the $X_i$, to ``fill space" by circumscribing an $n$-dimensional volume. A straightforward geometric interpretation of the vast majority of the other components remains to be discovered.

(\emph{Higher joint moments}). The covariance tensor (without qualification) is also the covariance tensor of trivial type $(l_i)$, where each $l_i=1$. If $n$ is small, the covariance tensor $t_{(l_i)}$ for non-trivial $(l_i)$ may be more interesting. In this case, the classical umbral calculus furnishes $GL(V)$ invariant functions of $t$ which may play a role similar to that of the Schur amplitudes. An explicit algorithm for such invariants is provided in \cite{grs}. Such invariants are a generalization of the $J$-invariants described in the special case of elements of $S^{4}\R^{2}$ or $S^{3}\R^{3}$.

\section{Algorithm}

In practice, information about the joint distribution of random variables $X_1,...,X_n$ is given by samples $(v^{j}_1,...,v^{j}_n)_{j=1..N}$, presumed to comprise an independent and identically distributed (i.i.d) sample from $(X_1,...,X_n)$ of size $N$.

\begin{definition}\label{samplecov}The \emph{sample covariance tensor} is
\begin{align*}
T=\widehat{\Cov}(X_1,...,X_n):=\sum_{j=1}^{N}(v^{j}_1-\bar{v}_1)\otimes\cdots\otimes (v^{j}_n-\bar{v}_n)\in V^{\otimes n}
\end{align*}
where $\bar{v}_i$ is the sample mean $\bar{v}_i=\tfrac{1}{N}\sum_{j=1}^{N}v^{j}_i$.

The \emph{sample Schur transform} or \emph{discrete Schur transform} is the tuple of components in the Schur-Weyl decomposition of $T$.
\end{definition}

This section provides the details of the numerical computation of the discrete Schur transform. The author's implementation is available at \url{http://github.com/schur-transform} .

\vspace{1pc}

\textbf{Precomputation}
\begin{enumerate}
\itemsep0em
\item{Fix a maximium size $n_{\text{max}}$ for the length of the input data series $\{(v^{j}_{1},...,v^{j}_{n})\}$, $j=1,...,N$, and dimension $k=\dim V$ for the data points, so that $v^{j}_{i}=(v^{j}_{i\alpha})_{\alpha=1,...,k}$ belongs to $\R^{k}$.}
\item{Compute the character tables of $S_n$ for each $n\leq n_{\text{max}}$.}
\item{Compute, for each $\sigma \in S_n$, the $k^{n}\times k^{n}$ matrix $P(\sigma)$ of permutation of tensor factors of $(\R^{k})^{\otimes  n}$ with respect to the basis $\{e_{\alpha_1}\otimes\cdots\otimes e_{\alpha_n}\}$, where $e_1,...,e_k$ is the standard basis of $\R^{k}$. Use lexicographical order of the index tuples $(\alpha_1,...,\alpha_n)$ for the basis ordering. For example, for the permutation $(12)$, the row corresponding to the input basis element $e_{\alpha_1}\otimes e_{\alpha_2}\otimes\cdots\otimes e_{\alpha_n}$ will have exactly 1 non-zero entry, the value 1 in the column corresponding to output basis element $e_{\alpha_2}\otimes e_{\alpha_1}\otimes\cdots\otimes e_{\alpha_n}$.}
\item{\label{perm}For each conjugacy class $c$ of $S_n$, compute the $k^n\times k^n$ matrix which is the sum given by:
\[S(c):= \sum_{\sigma \in c} P(\sigma) \]
}
\item{\label{projmatrix}For each character $\chi_\lambda$, i.e. each row of the character table of $S_n$, compute the $k^{n}\times k^{n}$ projection matrix which is the sum over conjugacy classes given by:
\[\pi(\lambda):= \sum_{\text{c}} \frac{\chi_\lambda(\text{identity})\chi_\lambda(c)}{n!}S(c) \]
}
\item{Verify that $\sum_\lambda \pi(\lambda)$ is the identity matrix. Since the entries of $n!\cdot \pi(\lambda)$ are integers, this equation should hold exactly.}
\end{enumerate}

\vspace{1pc}

\textbf{Main computation}
\begin{enumerate}
\itemsep0em
\item{A data series $v=(v^{j}_{i\alpha})$ is input, where $\alpha$ is the spatial dimension ranging from $1$ to $k$, $j$ is the sample dimension ranging from $1$ to $N$, and $i$ is the series dimension (e.g. time) ranging from $1$ to $n$.}
\item{Compute the means $\bar{v}_{i\alpha}:=\frac{1}{N}\sum_{j=1} v^{j}_{i\alpha}$, and replace $v^{j}_{i\alpha}$ with $v^{j}_{i\alpha}-\bar{v}_{i\alpha}$.}
\item{Compute the sample covariance tensor, a column vector of size $k^{n}$:
\[T:=\sum_{j=1}^{N}\thinspace\thinspace \sum_{\alpha_1,...,\alpha_n=1}^{k} \thinspace\thinspace\left(\prod_{i=1}^{n}  v^{j}_{i\alpha_i}\right)\quad e_{\alpha_1}\otimes\cdots\otimes e_{\alpha_n}\]}
\item{Compute the Schur transform, one tensor for each character $\chi_\lambda$, given by the matrix products:
\[T(\lambda)=\pi(\lambda)\cdot T\]}
\item{Verify that $T = \sum_\lambda T(\lambda)$.}
\item{Compute the Schur amplitude $|T(\lambda)|$, for each $\lambda$, as the square root of the sum of the squares of the numerical entries of $T(\lambda)$ (the Frobenius norm).}
\end{enumerate}

Note that in precomputation step \ref{projmatrix}, it is not necessary to use the inverse $\sigma^{-1}$ as in the formula of Theorem \ref{schurweyl}. In $S_n$, the conjugacy class $c$ of an element $\sigma$ is the permutation cycle type, which is the same for $\sigma$ and its inverse $\sigma^{-1}$. Thus $\chi_\lambda(\sigma)=\chi_\lambda(\sigma^{-1})=\chi_\lambda(c)$.

(\emph{Feasibility}). The algorithm presented above for the precomputation steps is not feasible for large $n$ because it requires iteration over the $n!$ elements of $S_n$. A more efficient algorithm could be obtained from a closed formula for the sum $S(c)$ over a given conjugacy class $c$ which does not require iteration over the members of $c$. However, even modestly large values of $n$ present a memory problem for the main computation steps, since the sample covariance tensor has $k^n$ numerical components and each projector $\pi(\lambda)$ has $k^{n}\times k^{n}$ components. For example, for $n=10$ and $k=3$, each projector uses over 10 GB of data storage. On the other hand, many of the Schur components $\pi_{\lambda}(V^{\otimes n})$ vanish entirely for certain dimensions of $V$ and the corresponding projectors do not need to be calculated (e.g. $\Lambda^{4}V=0$ if $\dim V \leq 3$).

\section{Applications}
\begin{figure}[h!]
  \caption{Violin plots of the 3-factor and 4-factor Schur content for a 6-timepoint time series of three-dimensional inhaling/exhaling lung CT scans, for 5 different subjects (numbered cases 1 through 5). The sequential/windowed 3-factor and 4-factor Schur content are also shown (in which only the consecutive 3-element and 4-element subsets of the time series are used). The data set is provided with $N=75$ expert-annotation landmarks, which were used as the samples.
  In all cases, the symmetric part, belonging to $S^{3}\R^{3}$ and $S^{4}\R^{3}$ respectively, labelled by the conjugacy class of the permutation $(123)$ in $S_{3}$ and $(1234)$ in $S_{4}$, dominates the Schur amplitude profile. This indicates that most of the points are moving very little, or are moving mainly radially with respect to the mean points. The volumetric part belonging to $\L^{3}\R^{3}$, labelled by the conjugacy class of the identity permutation $()$ in $S_3$, distinguishes cases 1 and 2, having much lower values, from cases 4 and 5, having much higher values, indicating that the trajectories in cases 1 and 2 are more planar with respect to the mean and the trajectories in cases 4 and 5 are more volumetric or 3-dimensional, deviating from a plane through the mean. On the other hand, the overall similarity between the profiles may indicate that the pattern appearing here is a ``signature" for breathing lung motion as opposed to, say, a head turn, a hand gesture, or a digestive tract motion.
  The CT data sets, described in \cite{castillo2009}, are available for download from \url{http://dir-lab.com}.}
  \centering
    \includegraphics[width=5.0in]{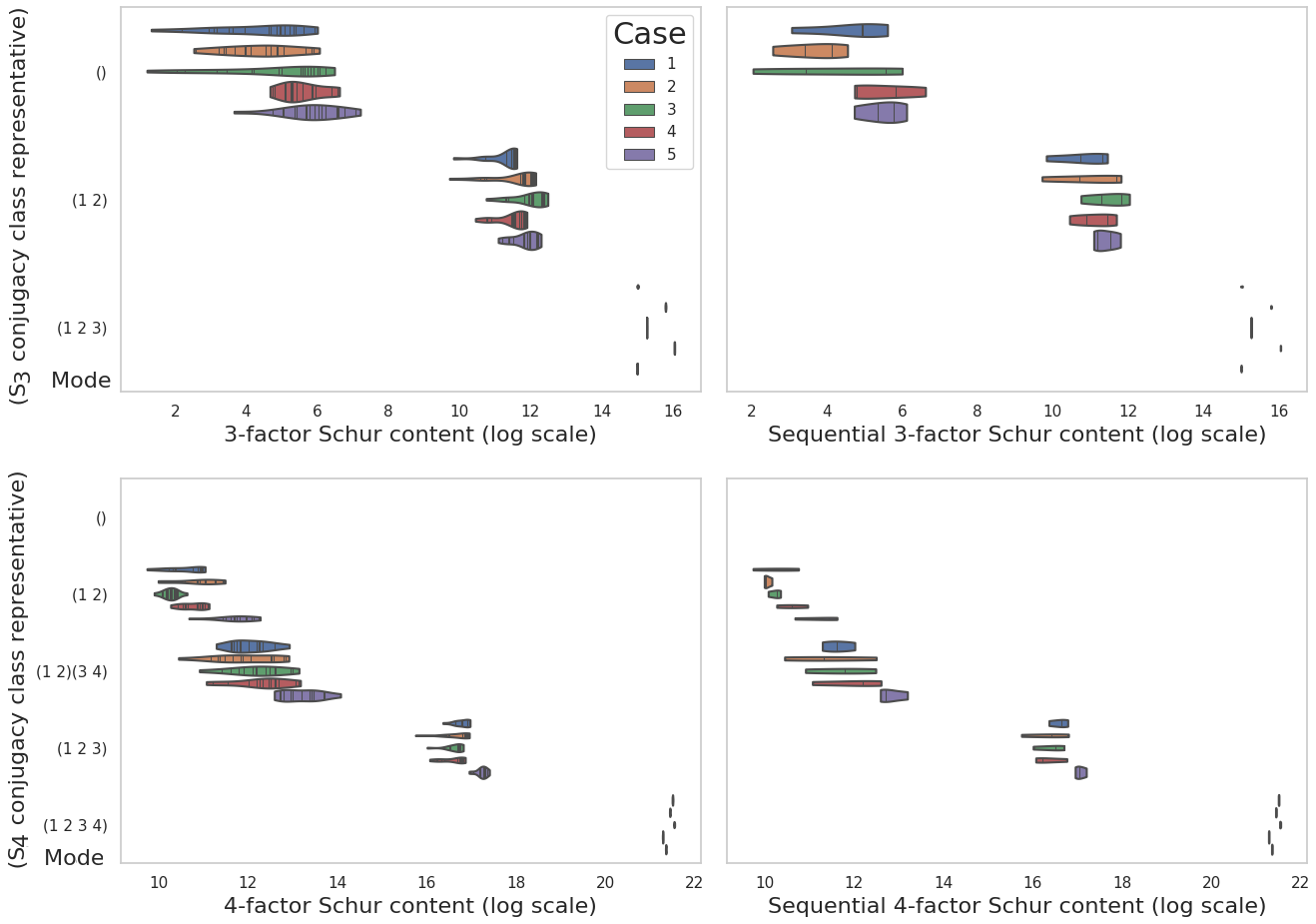}
\end{figure}

The Schur transform, or Schur content more generally, clearly has potential applications to morphometry, statistical shape analysis, fluid motion statistics, or body motion or gesture tracking.

(\emph{Exploratory classification and validation}). The $n$-factor Schur content summarizes the geometry of variation across a group of $m$ spatial variables ($m\geq n$). With increasing $n$, the resolution or complexity of the summary increases. Certain Schur components may be found to discriminate well between groups under different conditions, or certain patterns of Schur amplitudes may be found to characterize commonly occurring types of within-group variation.

(\emph{Classification rule}). Given a prior stratification of the $m$ matched objects of a given data set into classes, a possible classification rule for an additional $(m+1)$st object is as follows. For each class, evaluate both the $n$-factor Schur content $C$ and a modified $n$-factor Schur content $C'$ in which the $n$-fold subsets are replaced by the union of $(n-1)$-fold subsets with the additional sample. $C$ and $C'$ are both tuples of distributions, one for each $n$-partition type $\lambda$. Select the class which minimizes, for example, the $L^{1}$ or $L^{2}$ difference between the means of $C$ and the means of $C'$.

\bibliographystyle{plain}

\begin{thebibliography}{10}

\bibitem{berget}
Andrew Berget.
\newblock {\em Symmetries of Tensors}.
\newblock PhD thesis, University of Minnesota, 2009.

\bibitem{castillo2009}
Richard Castillo, Edward Castillo, Rudy Guerra, Valen Johnson, Travis McPhail,
  Amit Garg, and Thomas Guerrero.
\newblock A framework for evaluation of deformable image registration spatial
  accuracy using large landmark point sets.
\newblock {\em Physics in medicine and biology}, 54:1849--70, 04 2009.

\bibitem{fultonharris}
William Fulton and Joe Harris.
\newblock {\em Representation theory}, volume 129 of {\em Graduate Texts in
  Mathematics}.
\newblock Springer-Verlag, New York, 1991.
\newblock A first course, Readings in Mathematics.

\bibitem{grs}
Frank~D. Grosshans, Gian-Carlo Rota, and Joel~A. Stein.
\newblock {\em Invariant theory and superalgebras}, volume~69 of {\em CBMS
  Regional Conference Series in Mathematics}.
\newblock Published for the Conference Board of the Mathematical Sciences,
  Washington, DC; by the American Mathematical Society, Providence, RI, 1987.

\bibitem{tensormethodsstat}
Peter McCullagh.
\newblock {\em Tensor methods in statistics}.
\newblock Monographs on Statistics and Applied Probability. Chapman \& Hall,
  London, 1987.

\bibitem{schur1901}
I.~Schur.
\newblock {\em Ueber eine Klasse von Matrizen, die sich einer gegebenen Matrix
  zuordnen lassen}.
\newblock Dieterich in G{\"o}ttingen, 1901.

\bibitem{schur1927}
I.~Schur.
\newblock {\em Über die rationalen Darstellungen der allgemeinen linearen
  Gruppe}.
\newblock Sitzungsberichte Akad, 1927.

\bibitem{tuck}
Howard~G. Tucker.
\newblock {\em A graduate course in probability}.
\newblock Probability and Mathematical Statistics, Vol. 2. Academic Press,
  Inc., New York-London, 1967.

\bibitem{wass}
Larry Wasserman.
\newblock {\em All of statistics}.
\newblock Springer Texts in Statistics. Springer-Verlag, New York, 2004.
\newblock A concise course in statistical inference.

\bibitem{weyl1939classical}
H.~Weyl.
\newblock {\em The Classical Groups: Their Invariants and Representations}.
\newblock Princeton University Press, 1939.

\end{thebibliography}

\end{document}